\documentclass[12pt]{article}

\usepackage{amsmath,amsfonts}
\usepackage{paralist}
\usepackage{graphicx}  
\usepackage[colorlinks=true]{hyperref}
\usepackage{geometry} 
\usepackage[pagewise]{lineno}
\usepackage{overpic}
\usepackage{tikz}
\usepackage{cases}

\hypersetup{urlcolor=blue, citecolor=red}

\graphicspath{{./figs/}{./}}

\newcommand{\R}{\mathbb R}
\newcommand{\fmax}{f_{\textsc{max}}}
\newcommand{\fmaxhat}{\hat f_{\textsc{max}}}
\newcommand{\rhomax}{\tau}
\newcommand{\rhomaxmean}{\tau^{\textsc{ave}}}
\newcommand{\Dx}{\Delta x}
\newcommand{\Dt}{\Delta t}
\newcommand{\cfl}{\frac{\Delta t}{\Delta x}}

\newtheorem{rmk}{Remark}

\title{\bf \Large A Macroscopic Pedestrian Model with \\ Variable Maximal Density}

\author{
	Laura Bartoli\footnote{Dipartimento di Matematica e Fisica, Università Roma Tre, Rome, Italy.}, 
	Simone Cacace\footnote{Dipartimento di Matematica, Sapienza Università di Roma, Rome, Italy.}, 
	Emiliano Cristiani\footnote{Corresponding author; Istituto per le Applicazioni del Calcolo, Consiglio Nazionale delle Ricerche, Rome, Italy. \texttt{e.cristiani@iac.cnr.it}}, 
	Roberto Ferretti\footnote{Dipartimento di Matematica e Fisica, Università Roma Tre, Rome, Italy.}
}

\begin{document}
\maketitle

\begin{abstract}
\noindent In this paper we propose a novel macroscopic (fluid dynamics) model for describing pedestrian flow in low and high density regimes. 
The model is characterized by the fact that the maximal density reachable by the crowd -- usually a fixed model parameter -- is instead a state variable. 
To do that, the model couples a conservation law, devised as usual for tracking the evolution of the crowd density, with a Burgers-like PDE with a nonlocal term describing the evolution of the maximal density.
The variable maximal density is used here to describe the effects of the psychological/physical pushing forces which are  observed in crowds during competitive or emergency situations.\\
Specific attention is also dedicated to the fundamental diagram, i.e., the function which expresses the relationship between crowd density and flux. 
Although the model needs a well defined fundamental diagram as known input parameter, it is not evident \textit{a priori} which relationship between density and flux will be actually observed, due to the time-varying maximal density. An \textit{a posteriori} analysis shows that the observed fundamental diagram has an elongated ``tail'' in the congested region, thus resulting similar to the concave/concave fundamental diagram with a ``double hump'' observed in real crowds. \\
The main features of the model are investigated through 1D and 2D numerical simulations. The numerical code for the 1D simulation is freely available \href{https://gitlab.com/cristiani77/code_arxiv_2406.14649}{on this Gitlab repository}.
\end{abstract}


\section{Introduction}\label{sec:intro}
\paragraph{Context and motivations.} In this paper we deal with macroscopic (fluid dynamics) differential models for simulating pedestrian flow in normal and congested situations.
Generally speaking, computer simulations of crowd dynamics are of interest for getting insights about human group behavior and, more practically, for planning real mass events. 

Although macroscopic models are typically less used than microscopic differential models (like, e.g., social force models) or microscopic nondifferential models (like, e.g., cellular automata, optimal step models), they are still the only reasonable way to describe and forecast the behavior of very large-size crowd with thousands of individuals who frequently interact with their neighbors. 
The reason is that the computational cost of their implementation is \emph{independent of the crowd size}.

\paragraph{Relevant literature.}
The study of crowds is a multidisciplinary area which have attracted the interest of physicists, mathematicians, engineers, and psychologists.
Crowd modeling started from the pioneering papers \cite{hirai1975proc, okazaki1979TAIJp1, henderson1974TR} in the '70s.
Since then, all existing modeling approaches were investigated, spanning nanoscale, microscale (agent-based), mesoscale (kinetic), macroscale (fluid dynamic), and multiscale, either differential (based on ODEs and PDEs) or nondifferential (discrete choice, cellular automata, lattice gas), either continuous or discrete in time and space.
Also, models can be velocity based (first-order) or acceleration based (second-order), with local or nonlocal interactions, with metric or topological interactions. 
%
In addition to minor collision-avoidance actions (local navigation), models are also characterized by global path planning. This involves how pedestrians select their route to reach their destination, considering factors such as their familiarity with the environment, ability to predict movements, visibility conditions, and obstructions.
A number of review papers  
\cite{aghamohammadi2020TRB, 
	bellomo2011SR,
	bellomo2022M3AS,
	chen2018TR,
	corbetta2023AR,
	dong2020TITS,
	duives2013TRC,
	eftimie2018chapter, 
	li2019PhA,
	martinezgil2017CS,
	papadimitriou2009TRF,
	serena2023SMPT,
	yang2020GM}, 
meta-review papers
\cite{haghani2020SSp1, 
	haghani2020SSp2, 
	haghani2021PhA}
and books 
\cite{cristianibook, 
	rosinibook, 
	kachroobook, 
	maurybook} 
are available; we refer the reader to them for an introduction to the field. 
It is also worth to mention that crowd models often stem from, and share features with models developed in the context of vehicular traffic \cite{rosinibook, helbing2001RMP}.

More specifically, this paper has two sources of inspiration: 

1) The classical Hughes's model \cite{hughes2002TRB}, which is one of the most used first-order macroscopic models. It is based on the coupling of a conservation law, which describes the evolution of the density of people while they move towards a common target, and an Eikonal equation used to compute the minimum time path from any point of the domain to the target. 
It is also assumed that the minimum time path depends on the crowd's position itself, therefore the Eikonal equation must be repeatedly solved at any time. 
In this paper, we use a simplified approach, solving the Eikonal equation only once, assuming an empty domain (no people). This is done mainly for simplicity, but an extension to the fully coupled model is straightforward. 
The main similarity between the Hughes's model and the one proposed here is the presence of the fundamental diagram, i.e., the function which expresses the relationship between crowd's density (or speed) and flux. 
In both models a fundamental diagram must be specified and it is part of the input parameters. 

2) The interpersonal-distance model (IDM) recently introduced in \cite{cristiani2023PhA}, which is a nondifferential, continuous-in-space, discrete-in-time, microscopic (agent-based) model for pedestrian flow with cognitive heuristics. 
In the IDM there exist two state variables: the positions of the agents (as usual), and the minimum distance which agents accepts to keep between themselves and the other agents ahead. 
In other words, the agents want to keep a minimum distance from people they see (for comfort and collision-avoidance purposes), but they can also increase/decrease this threshold distance depending on the surrounding conditions. 
The rationale behind is that, in some competitive situations, it is preferable reducing comfort rather than risking to be overtaken by newly coming people, especially from behind.

\medskip

It is also useful to mention some papers about high density crowds, in both modeling and experimental literature.

From the experimental side, in \cite{helbing2007PRE, johansson2008ACS} the authors managed to compute the fundamental diagram until 10 ped/m$^2$. 
Remarkably, the fundamental diagram shows a concave/concave shape with a ``double hump'' (or a ``second peak''), namely the flux reaches two local maximum points. 
The double hump comes out because people are so densely packed that they start moving involuntarily, being pushed by the surrounding crowd (mass motion replaces individual motion) \cite{helbing2007PRE}.
See also \cite{jin2019TRC, lohner2017CD} for an experimental study that confirms the counter-intuitive increase of velocity at high densities.

From the modeling side, instead, we mention the macroscopic approach proposed in \cite{colombo2005M2AS} (see also \cite{chalons2007SISC}), where the mathematical properties of a fundamental diagram with double hump are investigated.
Velocity-based models with pushing forces are presented in \cite{kim2013SIGGRAPH, kim2015TVC} in the context of computer graphics, while social force models with pushing forces are presented in \cite{alrashed2020CD, helbing2000N, helbing2002PED, moussaid2011PNAS, yu2007PRE}. 
The paper \cite{moussaid2011PNAS} also adds cognitive heuristics.
A macroscopic model with crowd pressure force based on the Hughes's model is proposed in \cite{liang2021TRB}. 
A multiscale model for contact avoidance in high densities was investigated in \cite{narain2009SIGGRAPH}, while \cite{vantoll2021CeG} proposed another multiscale model based on SPH techniques. 
It also includes pushing forces and, remarkably, it is able to propagate material waves across the crowd triggered by pushing behavior.

%
%
%


In conclusion it is worth noting that the International Standard ISO 20414 was published in 2020 \cite{ISO20414}. This document outlines the procedures for verifying and validating evacuation models in the context of building fires, specifying 30 numerical tests to be conducted. Although this type of analysis is beyond the scope of the current paper, it would be beneficial to consider this document in future research.

\paragraph{Main contribution.} 
In this paper, we propose a first-order nonlocal macroscopic model for pedestrians, which reinterprets, in the 
differential context, some nice features of the IDM.
Since the IDM is nondifferential, there is no way to compute its many-particle limit using some consolidated analytical methods. 
Therefore, the model presented here is conceptualized from scratch, still it aims at preserving the foundational idea of the IDM, namely the presence of two state variables: the agents position, intended here in a macroscopic sense as the crowd \emph{density} $\rho$, and the minimum interpersonal distance allowed, intended here as the \emph{maximal reachable density} $\rhomax$. 
Note that the main novelty is in this point, since the maximal density reachable by the crowd -- usually a fixed model parameter -- is instead variable in space and time.
The dynamics of the maximal density $\tau$ is now driven by the evolution of an auxiliary variable, called $u$, which represents a sort of ``information'' propagating across the crowd in a domino effect. 
This information arises in competitive situations and translates
the urge for the crowd to ``compress'', if people do not want to be overtaken by newly coming people from behind, or to ``uncompress'', if there is no danger to lose the acquired priority within the crowd. 
The information can be either \emph{psychological}, if it comes from a simple observation of the crowdedness level, or \emph{physical}, if it comes from an actual pushing behavior.
In either case, we observe a ``material wave'' spanning the crowd, cfr.\ \cite{vantoll2021CeG}. 

As for the IDM, in the new model we observe two nice features: first, when a crowd reaches a steady state in front of a point of interest in a competitive scenario (like a crowd in a concert, when all people try to be as close as possible to the stage), we observe a self-emerging variable-in-space density, being higher near the point and lower far behind it, cf.\ \cite[Test 1a]{cristiani2023PhA}. 
This is a very realistic behavior which is hard to obtain with classical models, unless \emph{ad hoc} ingredients are added for this purpose, see, e.g.\ \cite{vonsivers2015TRB}.
Second, if we measure the relationship between density and flux \emph{ex post}, i.e., once the simulation is completed, we get a fundamental diagram with a long ``tail'' in the congestion region, which is very similar to that of IDM \cite[Test 3]{cristiani2023PhA} and to one experimentally measured by, e.g., Helbing et al.\ \cite{helbing2007PRE}. 
For these reasons, the proposed model is suitable to describe competitive situations at low and high densities.

Finally, the numerical tests will show that in the context of room evacuation with a bottleneck in front of the exit, our model seems to confirm the faster-is-faster effect \cite{adrian2020JRSI,haghani2019TRA}, rather than the celebrated faster-is-slower effect \cite{helbing2002PED}. 
This means that a rash evacuation, possibly with pushing behavior, actually \emph{reduces} the total evacuation time rather than increasing it.

\medskip 

The paper is organized as follows. In Sect.\ \ref{sec:model} we present and discuss the analytical model; Sect.\ \ref{sec:scheme} treats the numerical approximation, giving details for both the one- and the two-dimensional case. Lastly, Sect.\ \ref{sec:tests} presents and discusses some numerical tests aiming at highlighting the main features of the model.

\section{The model}\label{sec:model}
Let us denote by $x\in\R^2$ the space variable and by $t\in\R^+$ the time variable. We consider a crowd confined in a domain $\Omega\subset\R^2$, aiming at reaching a given target $\mathcal T\subset\partial \Omega$. We denote by $\rho=\rho(x,t)$ the crowd density at point $x$ and time $t$. 
It is natural to assume bounds for the density: we have $0\leq\rho(x,t)\leq\rhomax(x,t)$ for all $x$ and $t$, where $\tau$ represents the maximal density the crowd is allowed to reach. 
As already mentioned in the introduction, the model is characterized by the fact that, beside $\rho$, also $\rhomax$ has its own dynamics, coupled with that of $\rho$.

We also assume to know the shortest paths joining any point $x\in\Omega$ to the target $\mathcal T$. These paths can be easily found by solving once the stationary Eikonal equation
\begin{equation}\label{eikonalequation}
\left\{
\begin{array}{ll}
\|\nabla \phi(x)\|=1, & x\in\Omega \\
\phi(x)=0, & x\in\mathcal T \\
\phi(x)\to +\infty, & x\in\partial \Omega\backslash\mathcal T
\end{array}
\right.
\end{equation}
then integrating the vector field $w$ defined by
\begin{equation}\label{def:w}
w(x)=-\frac{\nabla\phi(x)}{\|\nabla\phi(x)\|},\quad x\in\Omega.
\end{equation}
Indeed, the solution of the equation 
\begin{equation}\label{ODEcontrollata}
	\left\{
	\begin{array}{l}
		\dot y(t)=w(y(t)) \\ 
		y(0)=x 
	\end{array}
\right.
\end{equation}
gives the shortest path from $x$ to $\mathcal T$ \cite[Sect.\ 8.2.3]{falconebook}.

We always assume that pedestrians look ahead while moving, where ``ahead'' means along the shortest path in the direction of the target, and ``behind'' means towards the opposite direction. 
This means that $w$ defines also the pedestrians orientation.

\begin{figure}[h!]
	\centering
\begin{tikzpicture}
	\draw[thick,->] (0,0) -- (9,0) node[anchor=north west] {$\rho$};
	\draw[thick,->] (0,0) -- (0,5) node[anchor=south east] {$f$};
	\node at (2,-0.3) {$\sigma$};
	\node at (-0.5,4) {$\fmax$};
	\node at (6.0,-0.3) {$\tau$};
	\node at (7.2,-0.3) {$\tau^*$};
	\node at (5,-0.3) {$\tau_*$};
	\node at (5.6,-0.3) {$<$};
	\node at (6.5,-0.3) {$<$};
	\draw[blue, thick] (0,0) -- (2,4);
	\draw[red, thick] (2,4) -- (7,0);
	\draw[red, thick] (2,4) -- (5,0);
	\draw[gray, thick, dashed] (2,4) -- (2,0);
	\draw[gray, thick, dashed] (2,4) -- (0,4);
\end{tikzpicture}
\caption{Fundamental diagram $\rho\mapsto f(\rho,\tau)$.}
\label{fig:FD}
\end{figure}
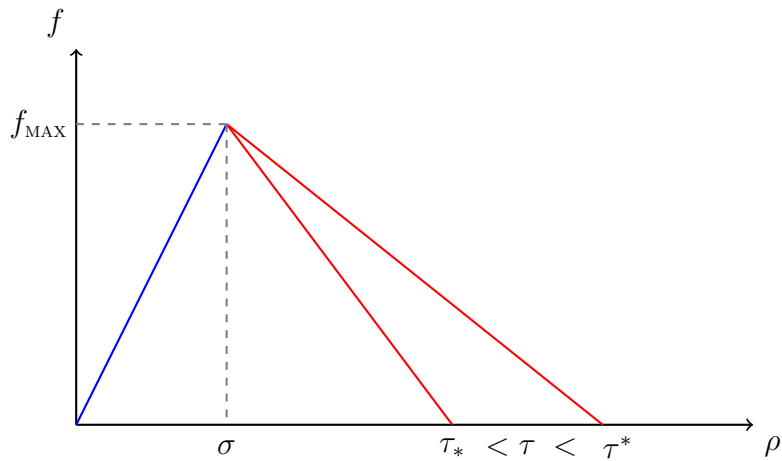

In order to present the model, we need to introduce an important ingredient, which is classical in macroscopic models: \emph{the fundamental diagram}. To this end, we assume that pedestrians at $x$ move in direction $w(x)$ with speed $s$. 
The speed $s$ is given as a function of the local density $\rho(x)$ and the maximal density $\rhomax(x)$, i.e., $s=s(\rho,\rhomax)$. 
Given the velocity field $s(\rho,\tau)w$, the flux function is 
$$
f=f(\rho,\rhomax,w)=\rho s(\rho,\rhomax)w
$$ 
and its modulus is said `fundamental diagram'. 
Several fundamental diagrams were considered in the literature \cite{vanumu2017ETRR}. 
In the following we choose a triangular-shaped function
\begin{equation}\label{def:flux}
	|f|=\rho s(\rho,\tau)=\left\{
	\begin{array}{ll}
		\frac{\fmax}{\sigma}\rho, & \rho\leq\sigma \\ [2mm]
		\frac{\rho-\rhomax}{\sigma-\rhomax} \fmax, & \rho>\sigma
	\end{array}
	\right.
\end{equation}
where $\fmax>0$ is the maximal flux and $\sigma$ is the critical density, see Fig.\ \ref{fig:FD}. Note that both $\fmax$ and $\sigma$ are fixed parameters, while $\rhomax$ is not.

\begin{rmk}
There is no theoretical reason to keep fixed $\fmax$ and $\sigma$. Indeed, it is reasonable to assume that they depend on $\tau$ as well. We made this choice just to highlight the role of the variable maximal density in the numerical tests.
\end{rmk}

Beside $f$, let us also define
\begin{equation}\label{def:g}
    g(u,w)=\frac12 u^2 w.
\end{equation} 

Let us introduce the last ingredients of the model: we assume that the dynamics is nonlocal, meaning that people are able to look far from them and take decisions on the basis of what they see. In particular, pedestrians evaluate the average maximal density ahead, defined by
\begin{equation}\label{def:rhomaxmean}
\rhomaxmean(x,t):=\frac{1}{|\mathcal S(x)|}\int_{\mathcal S(x)}\rhomax(z,t)dz,
\end{equation}
where $|\cdot |$ denotes the measure of the area and $\mathcal S(x):=\{z\in\Omega:|z-x|<\delta, \ w(x)\cdot (z-x)>0\}$ is the \emph{sensory region}, defined as the half ball of radius $\delta>0$ in front of each pedestrian.

We also introduce the following important threshold function
\begin{equation}\label{def:omega}
\theta[\rho,\tau](x,t):=\rho(x,t)-(\rhomaxmean(x,t)-\nu),
\end{equation}
where $\nu>0$ is a parameter (typically, small). 
The sign of $\theta$ is used to check if the current density $\rho$ has reached, within a tolerance $\nu$, the maximal average visible density $\rhomaxmean$.

We now have all the ingredients for presenting the complete model, which reads as
\begin{subnumcases}{\label{model}} 
	\partial_t\rho+\nabla\cdot (f(\rho,\tau,w))=0 \label{model-rho} \\ [1mm]
 \partial_t \rhomax=\gamma u \label{model-tau} \\ [0mm]
	\partial_t u+\nabla\cdot(g(u,w))=-\varepsilon u+
	\left\{
	\begin{array}{ll} 
		\alpha^+\Phi[\theta,w], & \theta\geq 0 \\ [1mm]
		\alpha^- \theta, & \theta<0
	\end{array}
	\right.
	\label{model-u} 
\end{subnumcases}
where $\alpha^+, \alpha^-,\gamma,\varepsilon>0$ are fixed parameters and $\Phi$ is a nonnegative function. 
For our purposes, it appears convenient to define
\begin{equation}\label{def:Phi}
\Phi[\theta,w]:= \max\{\theta -  \beta\nabla\theta\cdot w,0\}
\end{equation}
for some $\beta>0$. 
In addition, we have lower and upper bounds for $u$ and $\rhomax$, i.e., $u_*<u<u^*$ and $\rhomax_*<\rhomax<\rhomax^*$, with $\rhomax^*, \rhomax_*>0$ and $u^*>0$, $u_*<0$.
We also recall the natural constraint $\rho\leq\rhomax$.
Finally, the model is complemented with initial conditions $\rho(x,0)=\rho_0(x)$, $u(x,0)\equiv 0$, $\rhomax(x,0)=\rhomax_*$, and with boundary conditions, which will be discussed later on in the numerical tests.

Some comments are in order:
\begin{itemize}
    \item System \eqref{model} is fully coupled due to the dependence of $\theta$ (and hence $u$) on $\rho$ and $\tau$, see \eqref{def:omega}.
    \item The first equation \eqref{model-rho} is the classical conservation law describing the evolution of the density function along the vector field $w$ with speed $s(\rho,\tau)$. 
    \item The second equation \eqref{model-tau}, instead, provides the variation of the maximal density $\rhomax$ depending on the auxiliary function $u$.
    \item The third equation \eqref{model-u}, which is the main novelty of the model, is a Burgers-like equation with a source term and a reaction (annihilation) term.
    Here $u$ represents a sort of \textit{psychological} (in low density regimes) or \textit{physical} (in high density regimes) \textit{perturbation} which propagates as a wave through the crowd, more precisely along the shortest paths. 
    As in the classical Burgers equation, the speed of the propagation is given by height $u$ (positive or negative) of the perturbation itself. 
    Moreover, the perturbation locally arises and then travels towards the target, or in the opposite direction, according to the sign of $\theta$.
    In particular:
    \begin{itemize}
        \item if the sign of $\theta$ is positive, the density $\rho$ is approaching (from below) the density $\rhomaxmean$ (it is less-than-$\nu$ close to it), then a positive perturbation arises and starts traveling towards the target. While traveling, it transmits the information that some people are approaching from behind and, in practice, it increases the maximal density ahead. As a consequence, the people ahead pack themselves not to be overtaken and people behind can get a bit closer to the target, exactly as it happens in the IDM model  \cite{cristiani2023PhA};
        \item if the sign of $\theta$ is negative, the density $\rho$ is safely far from density $\rhomaxmean$, then a negative perturbation arises and starts traveling backward. While traveling, it transmits the information that there is more room ahead and, in practice, it decreases the maximal density. As a consequence, the people behind do not further pack themselves;
        \item the dissipation term $-\varepsilon u$ is added to dampen the perturbations, since it is unrealistic to assume that they propagate indefinitely across the crowd;
        \item Lastly, the most obscure term $\Phi$: let us first note that when $\Phi=0$, positive waves are no longer created. Now, assuming for the sake of clarity null boundary conditions for $u$ and constant boundary conditions for $(\rho,\tau)$, we see that, if both $\theta\geq 0$ (i.e., the density is close to the maximal one ahead) and $\Phi=0$, then $u$ is destined to vanish. As a consequence, $\tau$ will reach a steady state and $\rho$ will do the same (as in the classical LWR model for traffic flow).
        Roughly speaking, this corresponds to say that, if an equilibrium is reached, the function $\Phi$ determines the steady density profile. 
\end{itemize}
\end{itemize}

\begin{rmk}
	Following \cite{hughes2002TRB}, a Hughes-like version of the model can be easily created substituting \eqref{eikonalequation} with $s(\rho,\tau)\|\nabla \phi(x)\|=1$. Also, this would be essential to allow pedestrians to steer from the shortest path and walk along the minimum-time path instead, which can be time-dependent. 
    We have avoided to do that here to better highlight the role of the novel terms of the model.
\end{rmk}

\begin{rmk}
    One can note that if $\theta<0$ (i.e., the density is far from $\rhomaxmean$), negative waves for $u$ are continuously generated and tend to lower $u$ and then $\tau$ indefinitely. Actually, this does not happen only thanks to the lower bounds $u_*$ and $\tau_*$. Substituting $\alpha^\pm$ with $\rho\alpha^\pm$ we could prevent this to happen, at least when there is no crowd through which information can propagate. Again, we have avoided to do that to better highlight the role of the novel terms of the model.
\end{rmk}

\subsection{The model in 1D and its steady state}\label{sec:model1D}
In order to highlight some analytical properties of the model and to introduce the numerical approximation later on, it is convenient to consider the 1D version of the model.
Assuming to have $\Omega=[a,b]$ (a one-dimensional corridor) and $\mathcal T=\{b\}$ (therefore $w=1$), in 1D the model \eqref{model}, with the particular choice for $\Phi$ given in \eqref{def:Phi}, reads as
\begin{subnumcases}{\label{model1D}} 
	\partial_t\rho+\partial_x(\rho s(\rho,\tau))=0  \label{model1D-rho}\\ 
 \partial_t \rhomax=\gamma u \label{model1D-tau} \\
		\partial_t u+u\partial_x u=-\varepsilon u+
	\left\{
	\begin{array}{ll} 
		\alpha^+ \max\{\theta - \beta\partial_x\theta,0\}, & \theta\geq 0 \\
		\alpha^- \theta, & \theta<0.
	\end{array}
	\right. \label{model1D-u}
\end{subnumcases}

Let us now investigate a bit further the role of $\Phi$ in determining the steady state of the solutions in the congested regime ($\rho>\sigma$ in \eqref{def:flux} and also $\theta\geq 0$).
Consider a stationary regime characterized by $u\equiv 0$, which holds if the source term in \eqref{model1D-u} is identically zero. 
Let us ignore the trivial case $\theta-\beta\partial_x\theta<0$, and focus on the case $\theta-\beta\partial_x\theta=0$, so that
\begin{equation}\label{eq:theta_staz_1}
\theta(x) = c e^{x/\beta},\qquad c\in\R. 
\end{equation}
In order to have a stationary solution for $\rho$, we must require $f(\rho,\tau)$ to be constant, which amounts to a constraint between $\rho$ and $\tau$. 
Assuming $f(\rho,\tau)=d\>\fmax$, for some $d\in [0,1]$, Eq.\ \eqref{def:flux} gives
\begin{equation}\label{vincolo-rho-tau}
(\rho-\tau)=d(\sigma-\tau)
\qquad\Rightarrow\qquad
\rho = d\sigma + (1-d)\tau.
\end{equation}
Using now \eqref{eq:theta_staz_1} and \eqref{vincolo-rho-tau} into \eqref{def:omega}, we get
\begin{equation}\label{eq:staz}
c e^{x/\beta} = d\sigma + (1-d)\tau - \rhomaxmean + \nu.
\end{equation}
Then, a stationary solution, not necessarily unique, must satisfy \eqref{eq:staz} for a specific form of $\tau$ (which in turn implies the form of $\rhomaxmean$), and suitable values of the constants $c$ and $d$.

A first stationary solution may be obtained, for example, by setting $c=0$ (i.e., $\theta\equiv 0$), $d=0$ (i.e., $\rho=\tau$) and $\tau=ax+b$. 
In this case,
$$
\rhomaxmean = \frac 1 \delta \int_x^{x+\delta} (az+b) dz = ax + b + \frac{a}{2}\delta,
$$
which gives, once plugged into \eqref{eq:staz},
\begin{equation}\label{great_result_on_nu}
a=\frac{2\nu}{\delta}.
\end{equation}
Since $\nu,\delta>0$, we have $a>0$, therefore the density actually increases when getting closer to the target. 
Regarding the value of $b$, instead, it might be obtained by imposing a constraint on the total mass, i.e., on the integral of $\rho$.
Notably, this steady state seems to be the one found in the numerical solution, see, in particular, Test 2 in Section \ref{sec:tests} below.

\medskip 

A second form of stationary solution may be obtained assuming that $\tau(x)=a e^{x/\beta} + b$, which gives
\[
\rhomaxmean(x) = \frac 1 \delta \int_x^{x+\delta} \left(a e^{z/\beta} + b\right)dz = 
\frac{a\beta}{\delta}\left(e^{\delta/\beta}-1\right) e^{x/\beta} + b,
\]
so that using the above expression for $\tau$ and $\rhomaxmean$ in \eqref{eq:staz}, we obtain
\begin{equation}\label{eq:tau_exp}
c e^{x/\beta} =
d\sigma + (1-d)\left(a e^{x/\beta}+b\right) -
\frac{a\beta}{\delta}\left(e^{\delta/\beta}-1\right) e^{x/\beta}
-b + \nu.
\end{equation}
Then, \eqref{eq:tau_exp} splits into the system of conditions
\[
\begin{cases}
d(\sigma -b) + \nu = 0 \\
\displaystyle c=(1-d)a-
\frac{a\beta}{\delta}
\left(e^{\delta/\beta}-1\right).
\end{cases}
\]
Also in this second case a solution with $a>0$ (i.e., $\tau$ increasing) is admissible. 
Again, a condition on the total mass should be added.


%
%
%
%
%
%
%
%
\section{Numerical approximation}\label{sec:scheme}
\subsection{One-dimensional case}
We find it convenient to first expose the numerical approximation of equation \eqref{model} in 1D. 
This simplifies the notations and will allow to easily employ the splitting method \cite[Sect.\ 18.2]{levequebook} in the 2D case.

The two PDEs \eqref{model1D-rho}-\eqref{model1D-u} are both hyperbolic and can be solved using the Godunov scheme. 
Note that, even in 1D the flux function $f=\rho s(\rho,\tau)$ not only depends on the density $\rho$ but it depends also on space through the function $\tau$. 
Therefore, it is important to employ a generalized Godunov scheme which can deal with a space-dependent flux.

Let us introduce a numerical grid as usual, denoting by $\Dx$ the space step, by $\Dt$ the time step, by $j$ the index for space discretization, and by $n$ the index for time discretization.

For the first equation \eqref{model1D-rho} (concave flux) we have:
\begin{equation}\label{Godunov-rho}
	\rho_j^{n+1}=\rho_j^n-\cfl [\mathcal G^\rho(j,j+1,n)-\mathcal G^\rho(j-1,j,n)]
\end{equation}
where $\mathcal G^\rho$ is the numerical flux among two contiguous cells at time $n$, defined as
\begin{equation}\label{def:G^rho}
\mathcal G^\rho(j^\ell,j^r,n)=\min\{S^\rho(j^\ell,n),R^\rho(j^r,n)\}.
\end{equation}
Functions $S^\rho$ and $R^\rho$ are called \emph{sending} and \emph{receiving} functions, respectively, and are defined by
\begin{equation}\label{def:S-R^rho}
S^\rho(j,n)=
\left\{
\begin{array}{ll}
	f(\rho_j^n,\tau_j^n), & \rho_j^n\leq\sigma \\
	\fmax, & \rho_j^n>\sigma
\end{array}
\right.\quad , \quad
R^\rho(j,n)=
\left\{
\begin{array}{ll}
	\fmax, & \rho_j^n\leq\sigma \\
	f(\rho_j^n,\tau_j^n), & \rho_j^n>\sigma. 
\end{array}
\right.
\end{equation}

\medskip

For the second equation \eqref{model1D-u} (convex flux), let us denote by $\omega$ the right-hand side term and denote the flux by $g(u)=\frac12 u^2$. 
We have:
\begin{equation}\label{Godunov-u}
	u_j^{n+1}=u_j^n-\cfl [\mathcal G^u(j,j+1,n)-\mathcal G^u(j-1,j,n)]+\Dt\omega_j^n
\end{equation}
where $\mathcal G^u$ is the numerical flux among two contiguous cells, defined as
\begin{equation}\label{def:G^u}
\mathcal G^u(j^\ell,j^r,n)=\max\{S^u(j^\ell,n),R^u(j^r,n)\}
\end{equation}
with
\begin{equation}\label{def:S-R^u}
S^u(j,n)=
\left\{
\begin{array}{ll}
	0, & u_j^n\leq 0 \\
	g(u_j^n), & u_j^n>0
\end{array}
\right.\quad , \quad
R^u(j,n)=
\left\{
\begin{array}{ll}
	g(u_j^n), & u_j^n\leq 0 \\
	0, & u_j^n>0.  
\end{array}
\right.
\end{equation}

\medskip 

    While the scheme \eqref{Godunov-u} is nothing but a standard Godunov method for an equation with  source term, the scheme \eqref{Godunov-rho}, which corresponds to the Cell Transmission Method (CTM) \cite{daganzo1994TRB}, may be analysed in the framework of conservation laws with space-dependent fluxes. 
    Following \cite{zhang2005JCAM}, we can show that this definition of the scheme implies consistency and monotonicity.
	
    As for consistency, we recall from \cite{zhang2005JCAM} that the numerical flux $\mathcal G^\rho(j^\ell,j^r,n)$ is said to be {\em consistent} with the flux $f(\rho,x)=\rho s(\rho,\tau(x))$ if it is continuous with respect to $\rho_{j^\ell}$ and $\rho_{j^r}$, and satisfies the condition
	$$
	\mathcal G^\rho(j^\ell,j^r,n) = f(\rho_{j^\ell},x_{j^\ell}) = f(\rho_{j^r},x_{j^r})
	$$
	if $f(\rho_{j^\ell},x_{j^\ell}) = f(\rho_{j^r},x_{j^r})$, and in addition:
	\begin{enumerate}[(i)]
		\item
		$f_\rho(\rho_{j^\ell},x_{j^\ell}) f_\rho(\rho_{j^r},x_{j^r}) > 0$, or
		\item
		if $f_\rho(\rho_{j^\ell},x_{j^\ell})=0$, then $f_\rho(\rho_{j^r},x_{j^r}) \ge 0$;
		\item
		if $f_\rho(\rho_{j^r},x_{j^r}) = 0$, then $f_\rho(\rho_{j^\ell},x_{j^\ell}) \le 0$.
	\end{enumerate}
	In fact, if condition (i) holds, then, according to the definition of $S^\rho$ and $R^\rho$,
	$$
	\mathcal G^\rho(j^\ell,j^r,n)=\min\{S^\rho(j^\ell,n),R^\rho(j^r,n)\} = f(\rho_{j^\ell},x_{j^\ell}) = f(\rho_{j^r},x_{j^r}),
	$$
	since the min operator selects between $S^\rho$ and $R^\rho$ the one which equals the analytical flux. In the case of conditions (ii) or (iii), the form \eqref{def:flux} for the flux implies that $\rho_{j^\ell}=\rho_{j^r}=\sigma$ and therefore the previous consistency condition is trivially satisfied. Note that, in this latter case, we are implicitly assuming a regularization on $f$, so that the maximum of the flux would appear as a stationary point.
	
	On the other hand, monotonicity requires $\mathcal G^\rho$ to be nondecreasing with respect to $\rho_{j^\ell}$ and nonincreasing with respect to $\rho_{j^r}$. Since $S^\rho$ and $R^\rho$ are respectively nondecreasing and nonincreasing, the monotonicity condition results from the composition of monotone functions.

\subsection{Two-dimensional case}
The one-dimensional numerical approximation presented above, given in form of suitable combinations of sending and receiving functions, is already fit to deal with space-dependent fluxes.
Nevertheless, in 2D we face an additional difficulty: the fluxes 
$$
f=(f^x,f^y),\qquad g=(g^x,g^y)
$$ 
are no longer always positive in any of their components, but they can chance sign because of the dependence on $w=(w^x,w^y)$, see \eqref{def:w}.
This means that we can no longer assume that, e.g.\ along $x$ axis, the cell $j$ always sends mass to cell $j+1$, and that cell $j+1$ always receives mass from cell $j$. 
Instead, the direction of the flow now depends on the sign of $w^x$ and $w^y$. 

We proceed in this way: first, we apply the classical splitting procedure \cite[Sect.\ 18.2]{levequebook} in order to work separately in the $x$ and $y$ direction, thus unrolling a 2D problem in a series of 1D problems.

Starting with the equation \eqref{model-rho}, we have 
\begin{equation}
\partial_t \rho + \partial_x(\rho s w^x)+ \partial_y(\rho s w^y)=0.
\end{equation} 
Focusing on the motion along $x$ only, we start solving 
\begin{equation}\label{Godunov-rho*}
\partial_t \rho^* + \partial_x(\rho^* s w^x)=0,
\end{equation} 
for one time step, then we solve 
\begin{equation}
\partial_t \rho^{**} + \partial_y(\rho^{**} s w^y)=0
\end{equation} 
for another time step using $\rho^*$ as initial condition.

For \eqref{Godunov-rho*}, we define $\hat f^x=\rho s |w^x|$ so to have a positive flux (changing sign will be managed by duly selecting sending and receiving cells). 
Denoting by $k$ the grid index along $y$ axis, we have:
\begin{equation}\label{Godunov-rho-2D}
	\rho_{j,k}^{*}=\rho_{j,k}^n-\cfl [\mathcal H^\rho_{x,\textsc s}(j,k,n)-\mathcal H^\rho_{x,\textsc r}(j,k,n)]
\end{equation}
where 
\begin{equation}\label{def:HS}
\mathcal H^\rho_{x,\textsc s}(j,k,n):=
\left\{
\begin{array}{ll}
\mathcal G^\rho(j,j+1,k,n), & w^x\geq 0 \\ [2mm]
\mathcal G^\rho(j,j-1,k,n), & w^x<0.
\end{array}
\right. 
\end{equation}
and
\begin{equation}\label{def:HR}
\mathcal H^\rho_{x,\textsc r}(j,k,n):=
	\mathcal G^\rho(j-1,j,k,n)\mathbf{1}_{[0,+\infty)}(w^x_{j-1,k}) +
	\mathcal G^\rho(j+1,j,k,n)\mathbf{1}_{[0,+\infty)}(-w^x_{j+1,k}).
\end{equation}
where $\mathbf{1}$ is the indicator function.
In \eqref{def:HS} the mass leaving cell $(j,k)$ can be sent either to cell $(j+1,k)$ or $(j-1,k)$, depending on the sign of $w^x_{j,k}$. 
Analogously, in \eqref{def:HR} the mass enters cell $(j,k)$ from $(j-1,k)$, if $w^x_{j-1,k}>0$, and also from cell $(j+1,k)$, if $w^x_{j+1,k}<0$.

Similarly to \eqref{def:G^rho} and \eqref{def:S-R^rho}, we define
$$
\mathcal G^\rho(j^\ell,j^r,k,n)=
\min\{S^\rho(j^\ell,k,n),R^\rho(j^r,k,n)\},
$$
and
$$
S^\rho(j,k,n)=
\left\{
\begin{array}{ll}
	\hat f^x(\rho_{j,k}^n,\tau_{j,k}^n), & \rho_{j,k}^n\leq\sigma \\
	\fmaxhat, & \rho_{j,k}^n>\sigma
\end{array}
\right.\quad , \quad
R^\rho(j,k,n)=
\left\{
\begin{array}{ll}
	\fmaxhat, & \rho_{j,k}^n\leq\sigma \\
	\hat f^x(\rho_{j,k}^n,\tau_{j,k}^n), & \rho_{j,k}^n>\sigma .
\end{array}
\right.
$$
Analogous definitions are given for the second dimension $y$, i.e., $\hat f^y=\rho s |w^y|$, and $\mathcal H^\rho_{y,\textsc s}$, $\mathcal H^\rho_{y,\textsc r}$. 

All the construction must be also repeated for the equation \eqref{model-u} and the flux $g$, defining as before a positive flux $\hat g=(\hat g^x,\hat g^y)=\left(\frac12 u^2 |w^x|,\frac12 u^2|w^y|\right)$ and suitable $\mathcal H^u_{x,\textsc r}$, $\mathcal H^u_{x,\textsc s}$, $\mathcal H^u_{y,\textsc r}$, $\mathcal H^u_{y,\textsc s}$.
%
%
%
%
%
\section{Numerical tests}\label{sec:tests}
In this section we perform some 1D and 2D numerical tests to investigate the behavior of the model introduced above. 
Unless otherwise stated, the default set of parameters are those reported in Table \ref{tab:parameters}. Units of measure are `meters' and `seconds'. 
Note that, although we are mainly interested in the \emph{qualitative} behavior of the model, we have considered real-like parameters (e.g., the maximal speed of pedestrians is $\fmax/\sigma=1$ m/s).
\begin{table}[h!]
	\centering
	\caption{Default set of parameters used in the numerical tests}
	\begin{tabular}{|c|c||c|c|c|c|c|c|c|c|c|c|c|}
		\hline
		$\fmax$ & 
		$\sigma$ & 
		$\tau_*$ & 
		$\tau^*$ & 
		$u_*$ & 
		$u^*$ &
		$\varepsilon$ &
		$\alpha^+$ &
		$\alpha^-$ &
		$\beta$ &
		$\gamma$ &
		$\delta$ &
		$\nu$ \\ \hline\hline
		0.5 & 
		0.5 & 
		1 &
		5.5 &
		-1.5 &
		1 &
		0.1 &
		1 &
		0.1 &
		1 &
		0.01 &
		1 &
		0.1 \\
		\hline
	\end{tabular}
	\label{tab:parameters}
\end{table}
%

For the numerical approximation, we have considered $\Delta x=1$ ($\Delta x=0.5$ in Test 2, $\Delta x=2$ in Test 4b), and $\Delta t=\Delta x/2$ in order to fulfill the CFL condition.

\subsection{One-dimensional case}
In 1D, we consider a corridor where people move from left to right. 
At initial time $t=0$, a crowd with density $\rho=0.5$ is already present on the left side of a corridor of length 100, and continues to enter from the left side, with fixed density $\rho=0.5$, until time $t=150$. Instead, at the right boundary we set $\rho=0$ all the time.
A gate is located at $x_G=66$, and it is initially closed. 
After a while the crowd reaches the closed gate and starts queuing behind it. At time $t=400$ the gate opens and the crowd starts moving ahead, progressively reducing the queue until the end of the simulation.
Regarding boundary conditions for $u$ and $\tau$, we simply set them to 0 and $\tau_*$, respectively, on both sides of the domain. 
More important, the boundary conditions at gate: when it is closed, we set $\rho=\tau$ in the cell before it ($\tau$ and $u$ do not need specific Dirichlet conditions), and $(\rho=0, \tau=\tau_*, u=0)$ in all cells after it.

\begin{figure}[h!]
	\centering
	\includegraphics[width=0.49\textwidth]{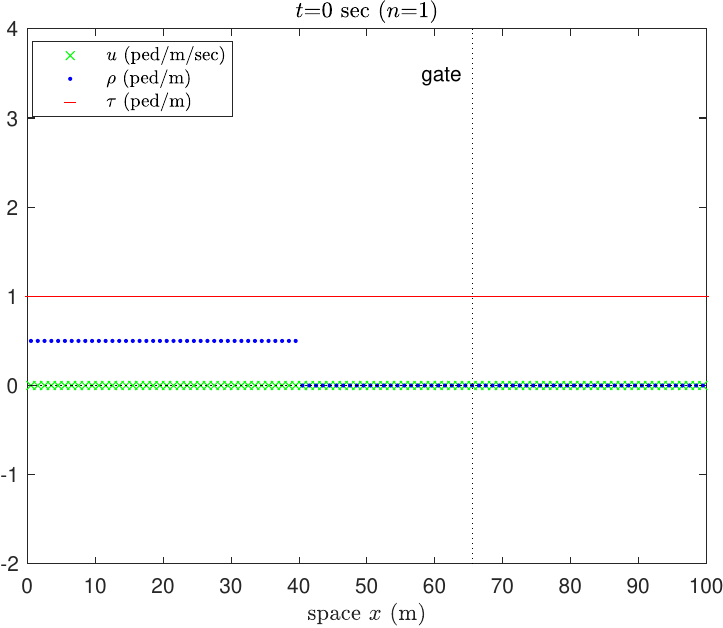}
	\includegraphics[width=0.49\textwidth]{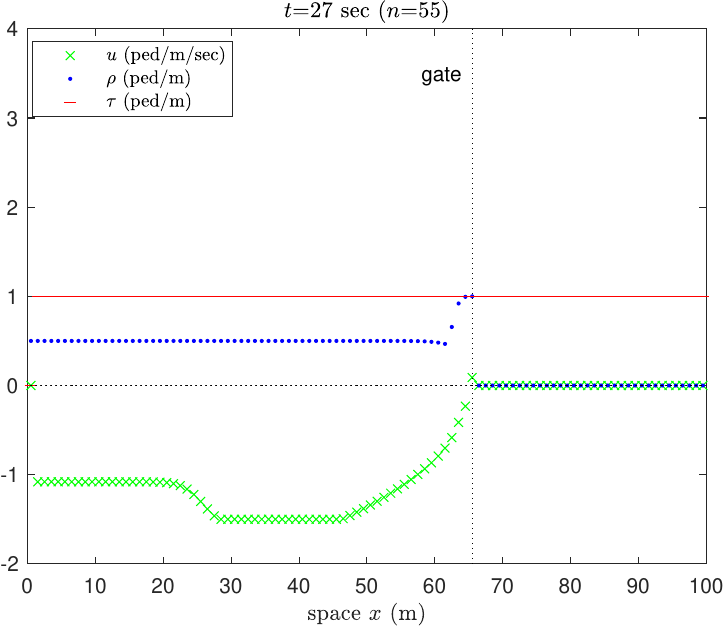} \\ [3mm]
	\includegraphics[width=0.49\textwidth]{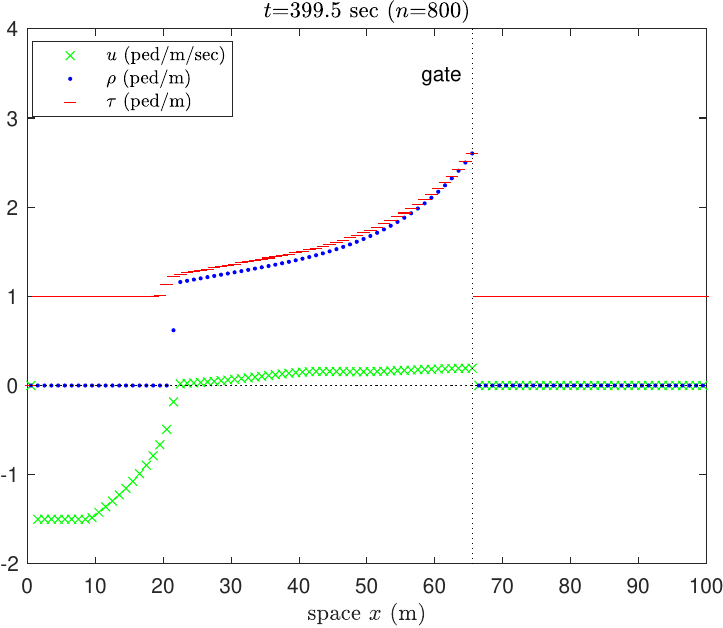}
	\includegraphics[width=0.49\textwidth]{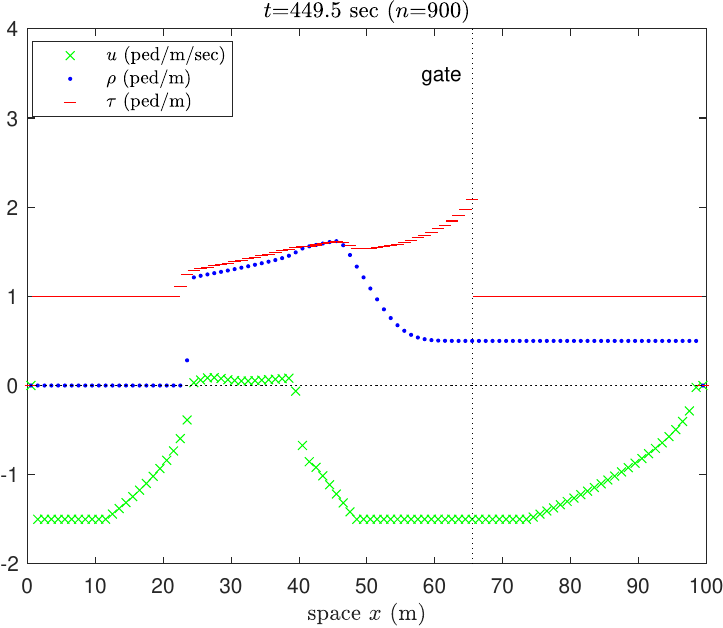}
	\caption{Test 1. Screenshots of the solutions $(\rho,\tau,u)$ at four time steps. From left to right, top to bottom: Initial condition, the queue behind the closed gate begins, the queue reinforces and back-propagates (with noncostant density), both the density and the maximal density drop after the gate opened. The numerical code for this test is freely available \href{https://gitlab.com/cristiani77/code_arxiv_2406.14649}{on this Gitlab repository}.
	}
	\label{fig:T1-screeshots}
\end{figure}

\paragraph{Test 1.}
Fig.\ \ref{fig:T1-screeshots} shows four screenshots of the simulation in the scenario described above. 
Here we have used $\varepsilon=0$ to better highlight the role of the PDE for $u$ (anyway the result with $\varepsilon=0.1$ is similar). 
We show the initial condition, the queue starts forming, the queue at its maximal length, and finally the crowd restarted after the opening of the gate.

It is interesting to see the behavior of $u$ and its effect first on $\tau$ and then on $\rho$. 
When $\rho$ increases and reaches the threshold $\rhomaxmean-\nu$, some positive waves are generated on $u$ (this happens around $n=50$). They start moving ahead increasing $\tau$ and then creating `more spaces' for the crowd. This makes $\rho$ increase as well. 
The most important feature of the model can be seen before the opening of the gate, when the queue has almost reached an equilibrium (see also Test 2 for the actual steady state configuration). 
In fact, the queue has a nonconstant density, being higher near the gate and lower far behind it. 
This is the same nice feature observed in \cite{cristiani2023PhA} in a microscopic setting.

When the gate opens, the crowd density $\rho$ starts decreasing, so negative waves appear on $u$ and they start propagating backwards. This makes the maximal density $\tau$ decrease as well and come back to the original state.

\begin{figure}[h!]
	\centering
	\includegraphics[width=\textwidth]{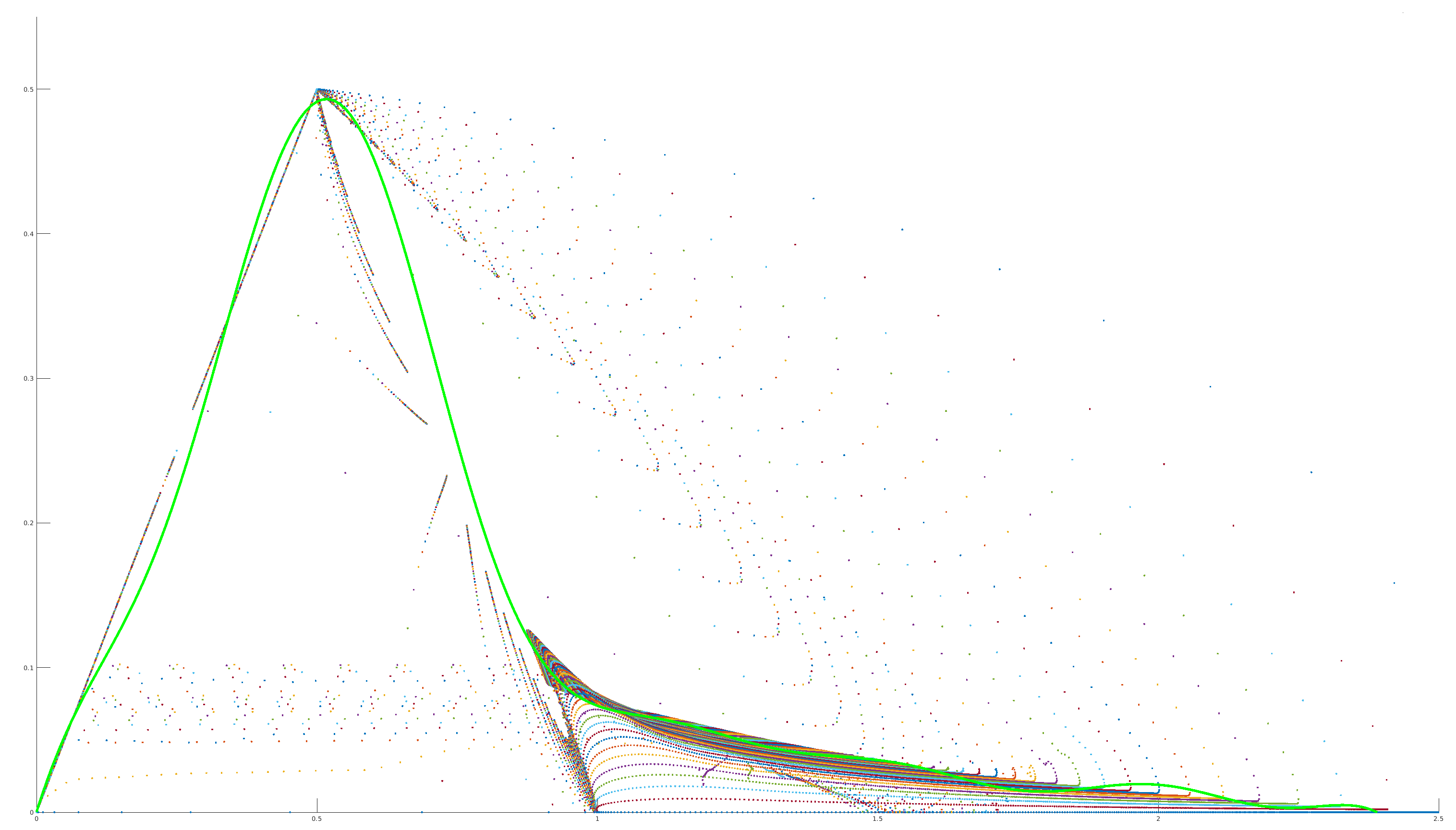}
	\caption{Test 1. Point-wise fundamental diagram computed \textit{ex post} by \eqref{FDexpost}. Each color is associated to a grid cell. The green line is the Matlab best fit obtained with sums of 5 sine's. The numerical code for this test is freely available \href{https://gitlab.com/cristiani77/code_arxiv_2406.14649}{on this Gitlab repository}.}
	\label{fig:T1-FD}
\end{figure}

In this simulation the density varies (increasing and decreasing) in a different ways in all grid cells. Therefore, it is interesting to recover, \emph{ex post}, the shape of the fundamental diagram plotting all pairs density-flux 
\begin{equation}\label{FDexpost}
	(\rho, G^\rho(j,j+1,n))
\end{equation}
at any time step $n$ and every cell $j$. 
This analysis must be regarded as a sort of `experimental' measurement performed on the simulation results, exactly in the same spirit of the evaluation done in \cite{cristiani2023PhA} where authors obtain a double-hump fundamental diagram. Moreover, the result can be compared with actual experimental measurements obtained, e.g., in \cite{helbing2007PRE}. 
Fig.\ \ref{fig:T1-FD} shows the result.

Since this is a macroscopic model based on the specific fundamental diagram \eqref{def:flux}, obviously we expect to recover that shape in the analysis. Indeed, the shape of the triangular function is clearly visible. Nevertheless, it is interesting to see that the triangle has a sort of `tail' on the right side which comes from the fact that the maximal density can vary from $\tau_*$ to $\tau^*$.   

The numerical code for this test is freely available \href{https://gitlab.com/cristiani77/code_arxiv_2406.14649}{on this Gitlab repository}.

\begin{figure}[h!]
	\centering
	\includegraphics[width=0.49\textwidth]{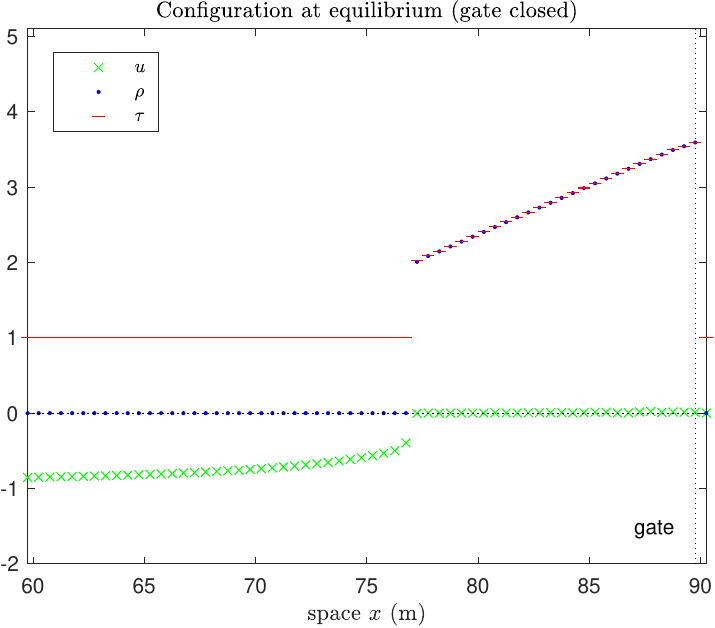}
	\includegraphics[width=0.49\textwidth]{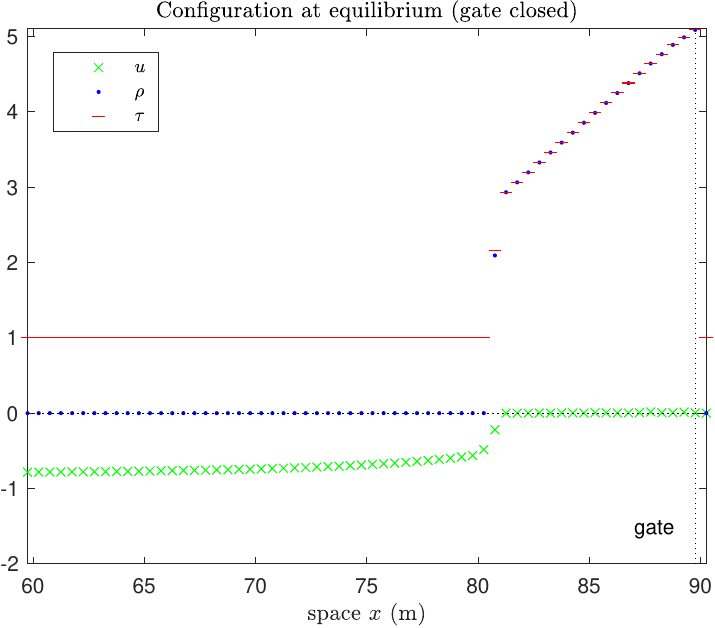} 
	\includegraphics[width=0.49\textwidth]{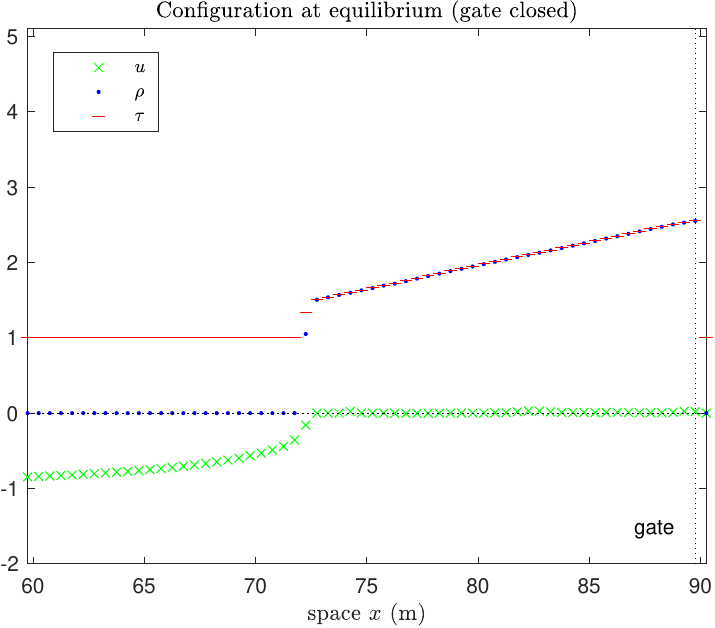} 
	\includegraphics[width=0.49\textwidth]{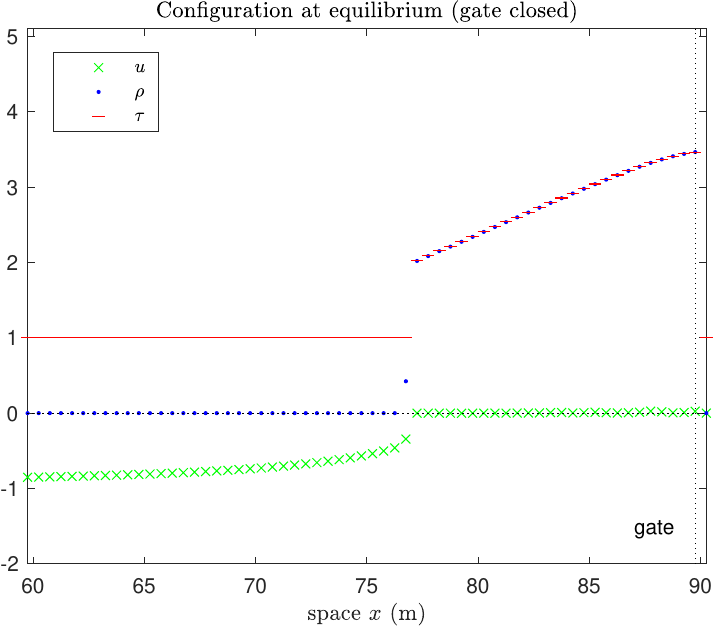}
	\caption{Test 2. Screenshots of the solutions $(\rho,u,\tau)$. From left to right, top to bottom: default set of parameters (Table \ref{tab:parameters}), $\nu$ increased to 0.2, $\nu$ decreased to 0.05, $\beta$ increased to 2.}
	\label{fig:T2-sensitivity}
\end{figure}

\paragraph{Test 2.}
In this test we keep the gate closed all the time, allowing the queue to reach a steady state.
Fig.\ \ref{fig:T2-sensitivity} shows the final configuration obtained by using the default set of parameters (Table \ref{tab:parameters}) as well as varying some of them. Variations are useful for a sensitivity analysis which highlights the role of the different parameters.
Notably, it seems that we get the linear steady state for $\tau(x)$ analytically studied in Sect.\ \ref{sec:model1D}.

One can note that $\nu$ is the most effective parameter to adjust the slope of the crowd density, see \eqref{great_result_on_nu}. $\delta$ has also a similar effect, but likely introduces instabilities. $\beta$, instead, changes the slope of the queue at the right-end.

Note that variations on $\alpha^\pm$ and $\gamma$ are not shown since they mainly rule the temporal scale (if not create instabilities).
\clearpage
\subsection{Two-dimensional case}
In 2D we consider a simpler scenario where a crowd has to evacuate a square domain of size 100 $\times$ 100 as fast as possible, through one or two exits located on the right side.
Exits are large as one numerical cell, i.e., $\Delta x$.

\begin{figure}[h!]
	\centering
	\includegraphics[width=0.32\textwidth]{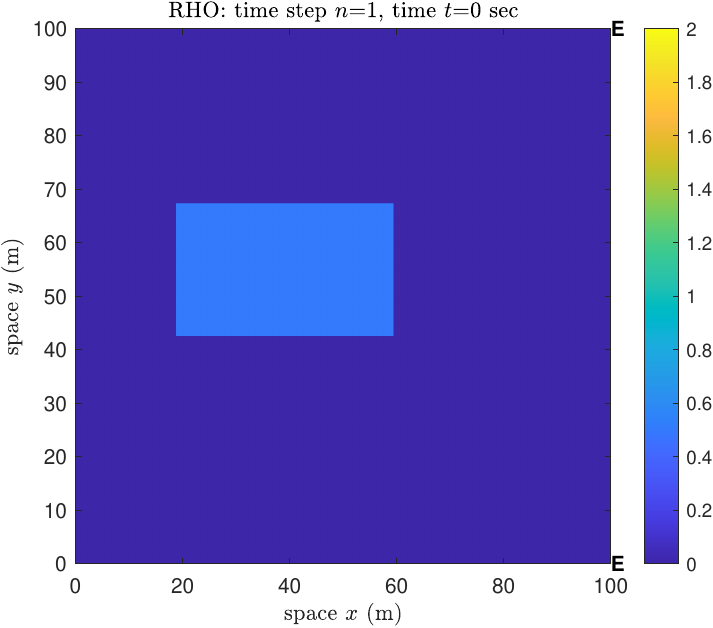} 
        \includegraphics[width=0.32\textwidth]{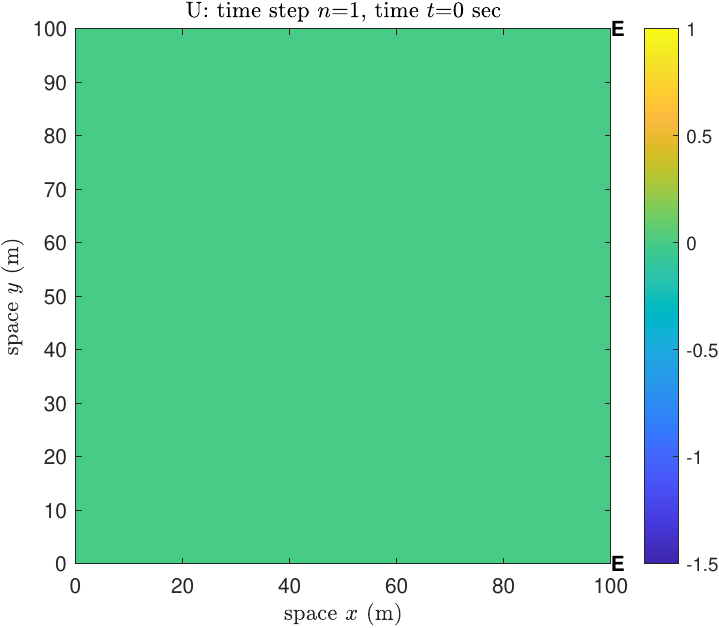} 
	\includegraphics[width=0.32\textwidth]{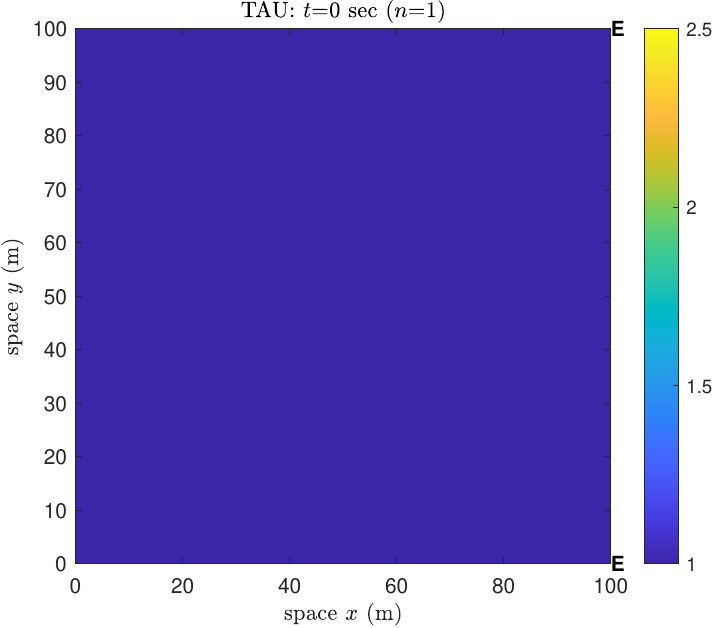} 
        \includegraphics[width=0.32\textwidth]{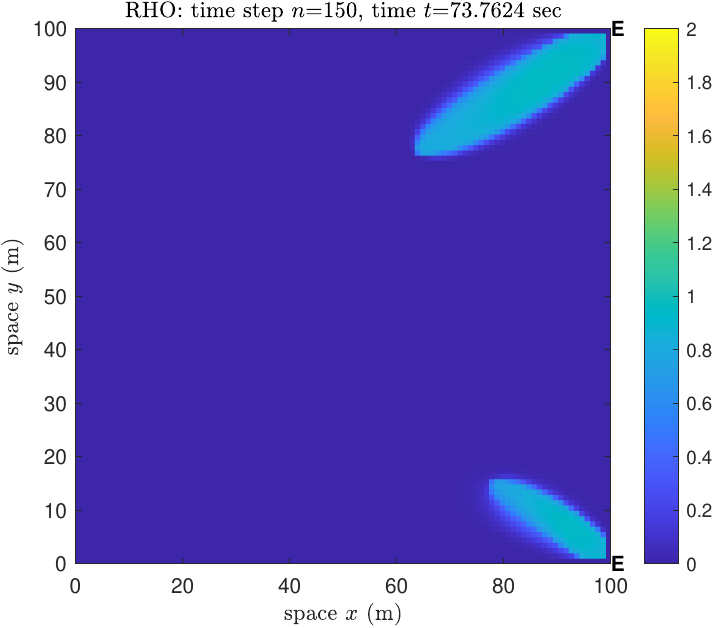} 
        \includegraphics[width=0.32\textwidth]{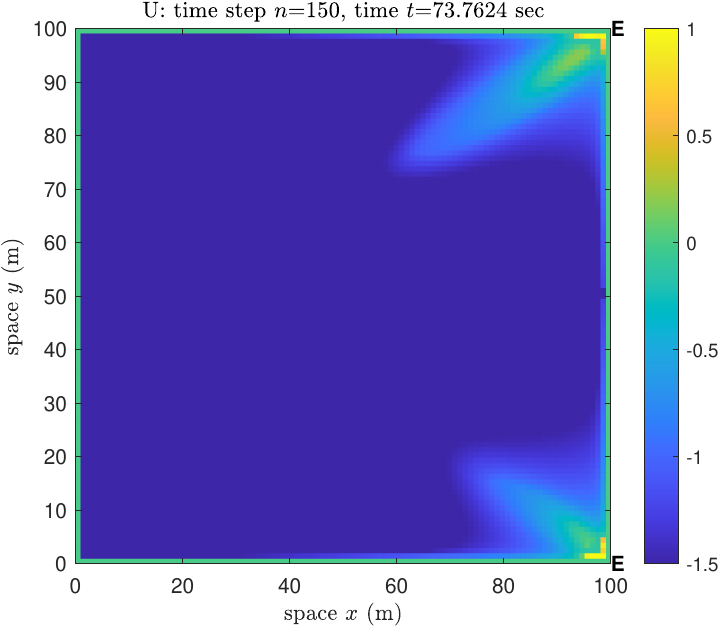} 
	\includegraphics[width=0.32\textwidth]{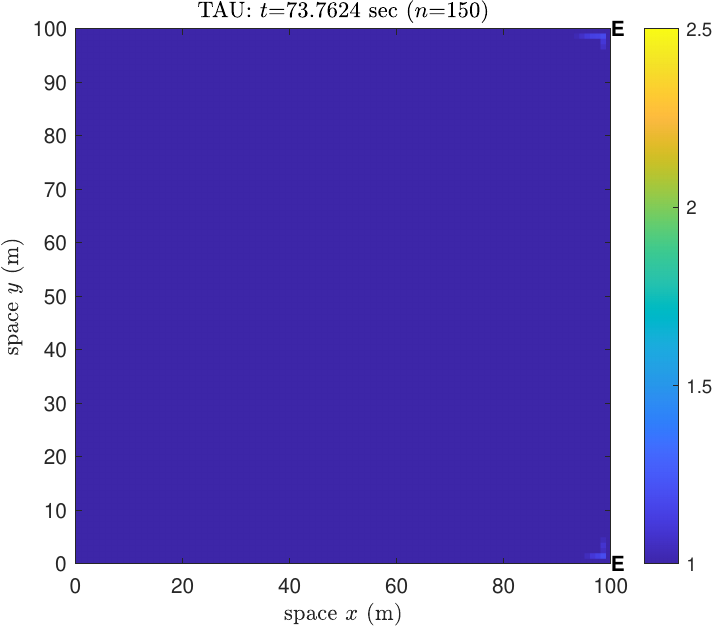} 
        \includegraphics[width=0.32\textwidth]{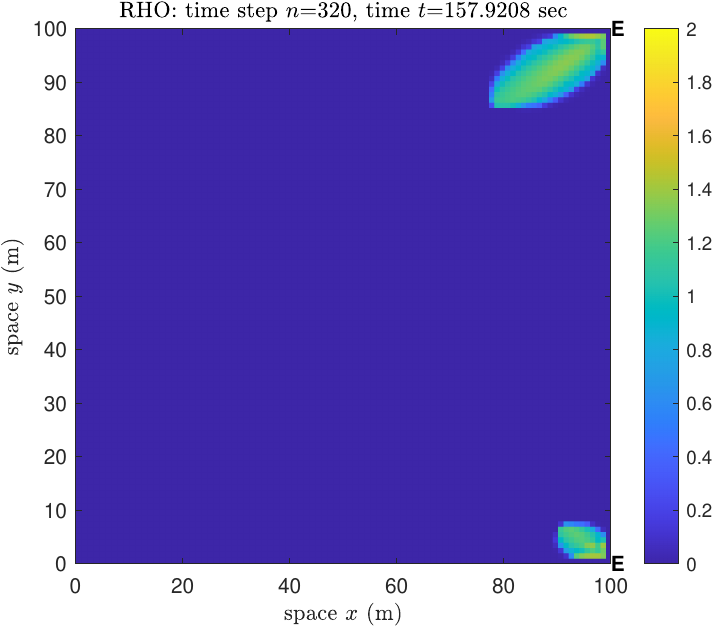} 
        \includegraphics[width=0.32\textwidth]{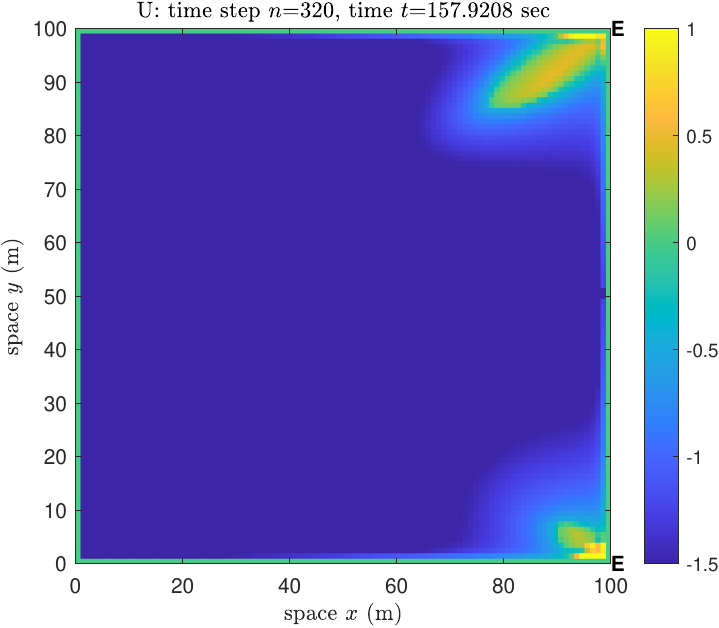} 
	\includegraphics[width=0.32\textwidth]{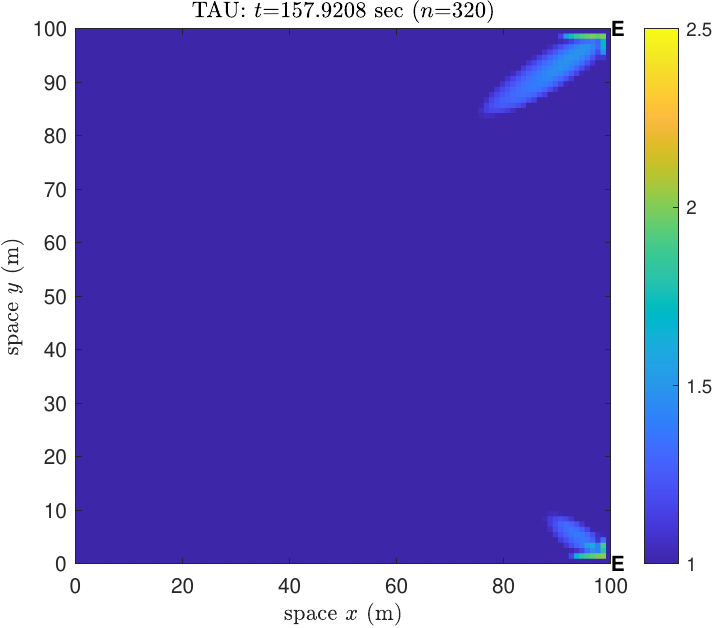} 
	\caption{Test 3. Screenshots of the solutions $\rho$ (first column), $u$ (second column), and $\tau$ (third column), at $t=0$, 75, 160.}
	\label{fig:T3-doubleexit}
\end{figure}

\paragraph{Test 3.}
In this test the two exits are located at the top- and bottom-right corners.
At initial time $t=0$ the crowd is in $[20,60]\times[44,68]$ and has constant density $\rho=0.5$. 
Fig.\ \ref{fig:T3-doubleexit} shows three screenshots.

The crowd immediately splits horizontally (in two unequal parts) in such a way that each person heads toward the closest exit.
Similarly to the 1D case, when people approach the exit, the density $\rho$ crosses the threshold $\rhomaxmean-\nu$, and positive waves arise on $u$. Waves travel along the shortest paths, hence directed toward the exits. The waves make the maximal density $\tau$ increase and then the density $\rho$ increases as well.

\begin{figure}[h!]
	\centering
	\includegraphics[width=0.32\textwidth]{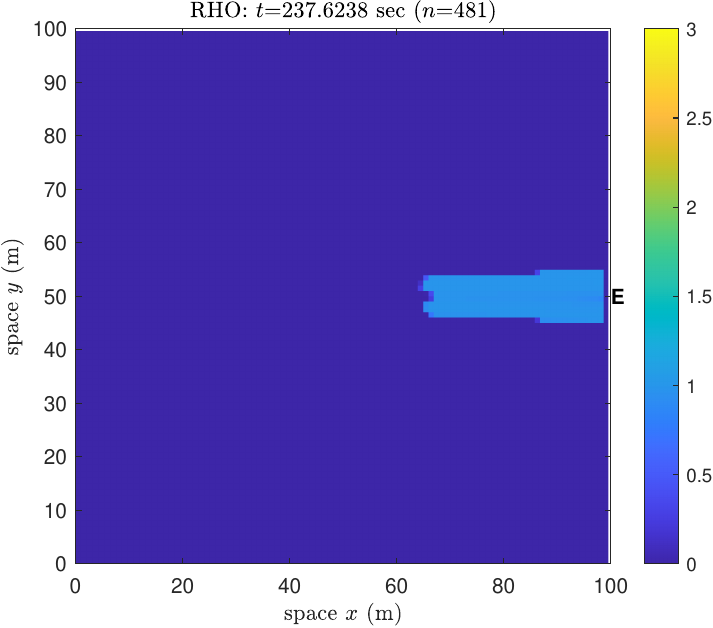} 
        \includegraphics[width=0.32\textwidth]{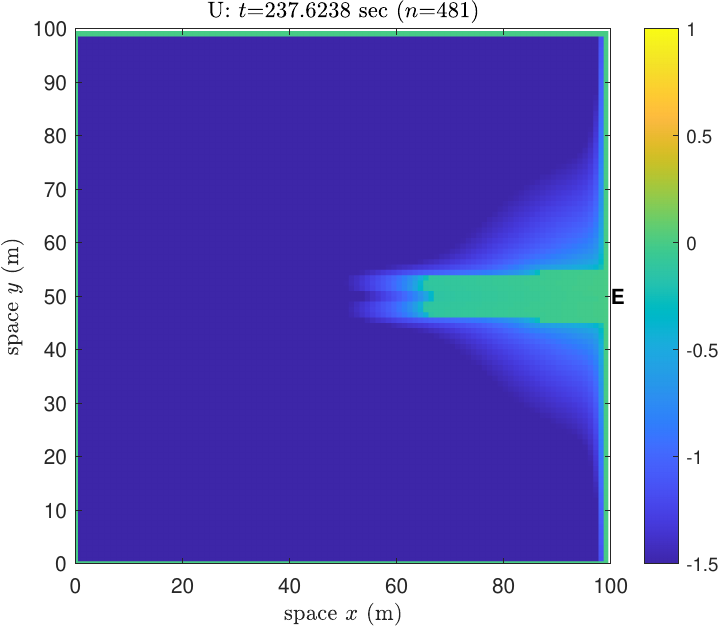}
        \includegraphics[width=0.32\textwidth]{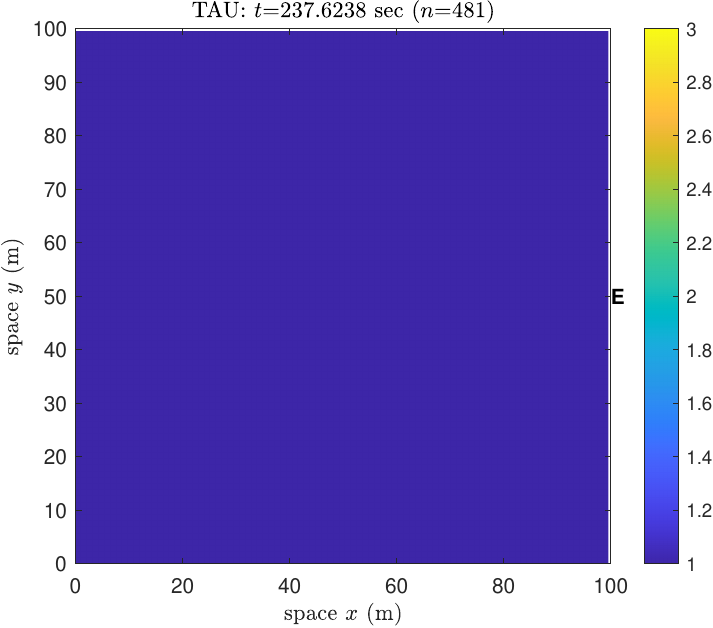}
	\includegraphics[width=0.32\textwidth]{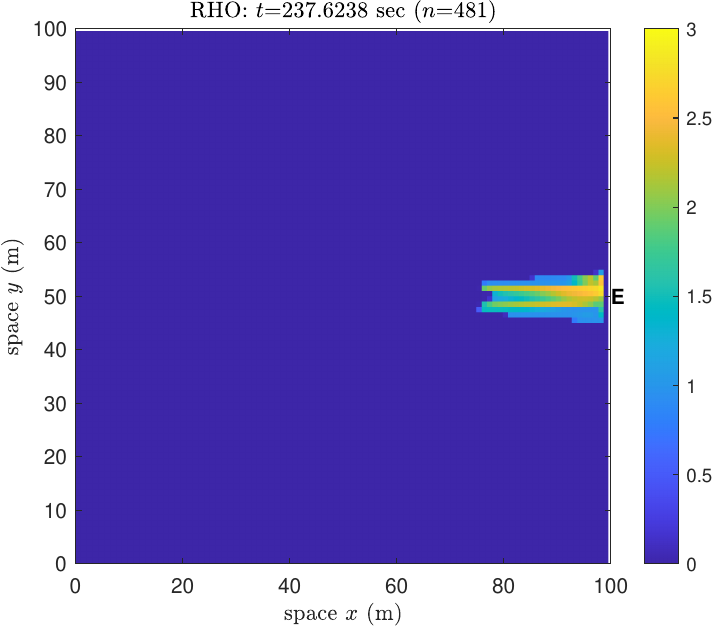} 
        \includegraphics[width=0.32\textwidth]{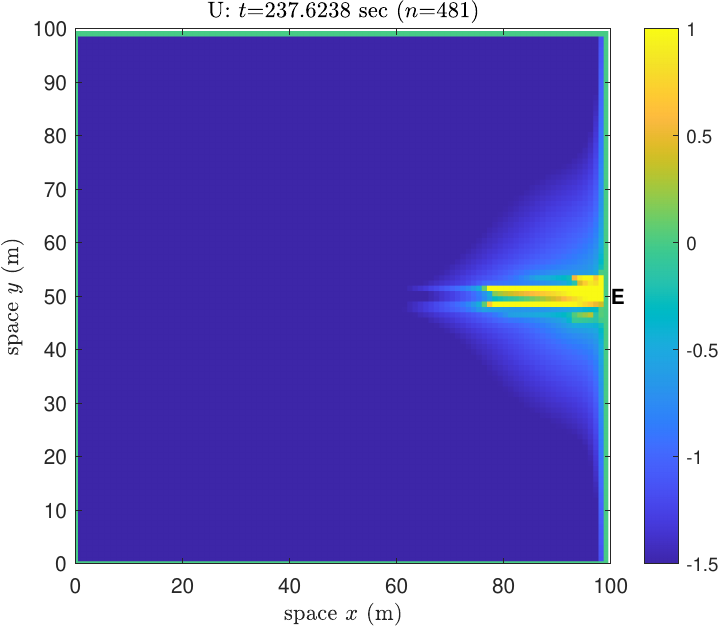}
        \includegraphics[width=0.32\textwidth]{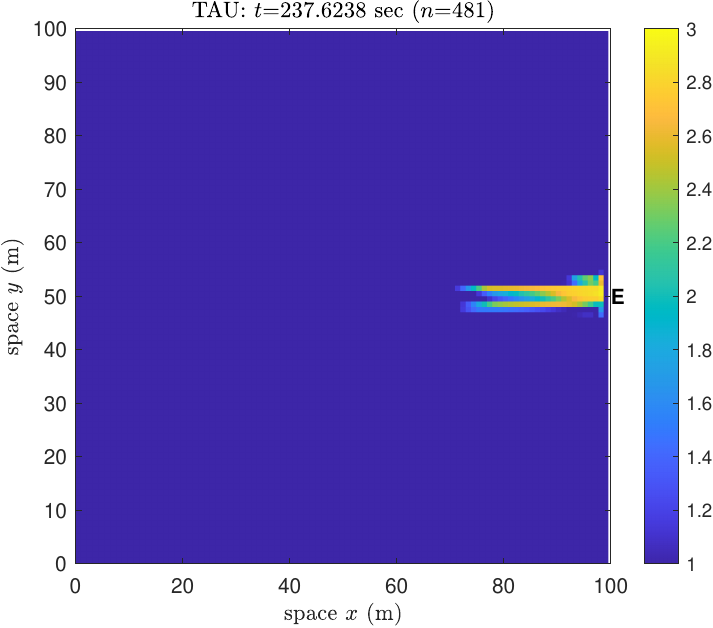}
	\caption{Test 4a. Screenshots of the solutions $\rho$ (first column), $u$ (second column), and $\tau$ (third column), at $t=240$. $\alpha^+=0$ (top), $\alpha^+=1$ (bottom).}
	\label{fig:T4a-singleexit_screenshots}
\end{figure}

\paragraph{Test 4.}
In this test we consider only one exit in the middle of the right side of the domain and slightly different initial conditions for $\rho$. 
More important, we consider two kind of obstructions for the crowd: 
\\
\\
\indent $\bullet$ Test 4a: we assume that the outgoing flux at the exit is reduced by a factor 1/2.
Fig.\ \ref{fig:T4a-singleexit_screenshots} shows two screenshots of the solutions captured at the same time $t=240$ and obtained with $\alpha^+=0$ and $\alpha^+=1$, respectively. 
Both the evacuation rate and the total evacuation time are identical in the two cases, but one can see that the queue's shape differs. 
If $\alpha^+=0$, no perturbation arises on $u$, then the maximal density never changes. 
If $\alpha^+=1$, instead, the usual mechanism $\rho$ increases $\Rightarrow$ $u$ increases and moves $\Rightarrow$ $\tau$ increases $\Rightarrow$ $\rho$ increases, is triggered.  

\medskip 

\begin{figure}[h!]
	\centering
	\includegraphics[width=0.5\textwidth]{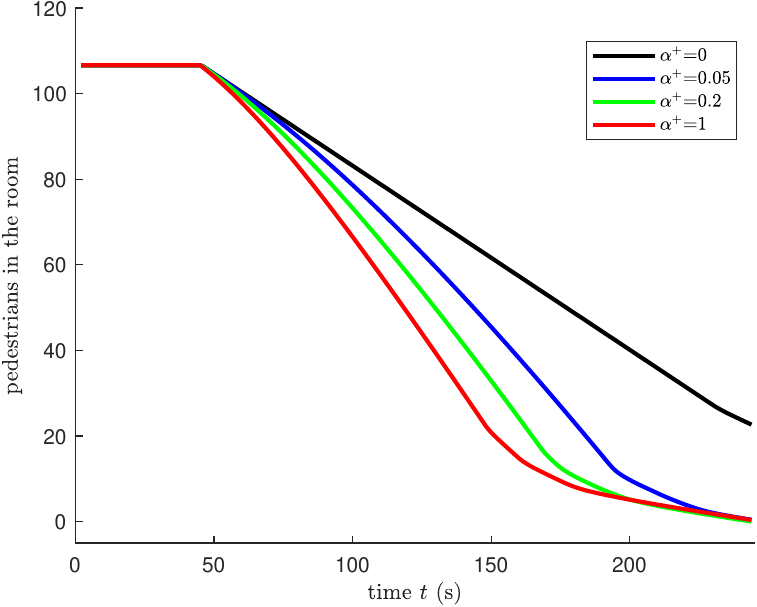}
	\caption{Test 4b. Pedestrian mass $\iint \rho \ dx dy$ inside the room as a function of time, for $\alpha^+=0,0.05,0.2,1$.}
	\label{fig:T4b-singleexit_mass}
\end{figure}

$\bullet$ Test 4b: we consider an obstacle $\Delta x$ wide just in front of the exit, obtained by imposing a Dirichlet condition $\rho=\rhomax_*-0.1=0.9$ at the cell in front of the exit.
Conversely to Test 4a, here the queue's shape remains basically the same but the evacuation rate does depend on $\alpha^+$. 
Fig.\ \ref{fig:T4b-singleexit_mass} shows four evacuation rates obtained for $\alpha^+=0,\ 0.05,\ 0.2,\ 1$. 

One can see that a larger $\alpha^+$ speeds up the evacuation. 
As already anticipated at the end of 
Sect.\ \ref{sec:intro}, this result is in line with the so-called `faster-is-faster' effect, which basically states that a rush evacuation, possibly with pushing behavior, actually decreases, instead of increasing, the evacuation time; see the recent papers \cite{adrian2020JRSI,haghani2019TRA}, which states that the evacuation time is actually reduced by a greater `pressure' along the crowd.

\small 
\paragraph{Funding.}
E.C.\ would like to thank the Italian Ministry of University and Research (MUR) to support this research with funds coming from PRIN Project 2022 PNRR (No. 2022XJ9SX, entitled ``Heterogeneity on the road - Modeling, analysis, control'').

E.C.\ would like to thank the Italian Ministry of University and Research (MUR) to support this research with funds coming from PRIN Project 2022 (No. 2022238YY5, entitled ``Optimal control problems: analysis, approximation and applications'').

This study was carried out within the Spoke 7 of the MOST -- Sustainable Mobility National Research Center and received funding from the European Union Next-Generation EU (PIANO NAZIONALE DI RIPRESA E RESILIENZA (PNRR) – MISSIONE 4 COMPONENTE 2, INVESTIMENTO 1.4 – D.D. 1033 17/06/2022, CN00000023). This manuscript reflects only the authors' views and opinions. Neither the European Union nor the European Commission can be considered responsible for them. 

E.C.\, R.F., and S.C. are members of the Gruppo Nazionale Calcolo Scientifico - Istituto Nazionale di Alta Matematica (GNCS-INdAM).


\bibliographystyle{unsrt}
\bibliography{biblio}
\end{document}